\documentclass{czjphys}         
\usepackage{amsmath,amssymb,amsthm}

\newtheorem{thm}{Theorem}[section]
\newtheorem{prop}[thm]{Proposition}

\newtheorem{cor}[thm]{Corollary}
\theoremstyle{definition}
\newtheorem{defn}[thm]{Definition}
\newtheorem{eg}[thm]{Example}
\newtheorem{rem}[thm]{Remark}

\newcommand{\Ah}{A_\hb}
\newcommand{\ad}{{\mathrm{ad}}}
\newcommand{\al}{\alpha}
\renewcommand{\be}{\beta}
\newcommand{\clal}{{\ol{\al}}}
\newcommand{\clbe}[1][]{{\ol{\be}_{#1}}}
\newcommand{\De}{\Delta}

\newcommand{\Ga}{\Gamma_\k}
\newcommand{\g}{\mathfrak g}

\newcommand{\h}{\mathfrak h}
\newcommand{\hb}{h}
\newcommand{\hd}{{\mathfrak{h}}^*}
\renewcommand{\k}{\mathfrak{k}}
\newcommand{\ka}{\varkappa}
\newcommand{\Ker}{{\mathrm{Ker}}}
\newcommand{\kk}{\mathbb{C}}
\newcommand{\la}{\lambda}
\newcommand{\mh}{\mu_\hb}
\newcommand{\Omg}{\Omega}
\newcommand{\ol}{\overline}
\newcommand{\ph}{\varphi}

\newcommand{\Tor}{{\mathrm{Tor}}}
\newcommand{\U}{{\mathrm{U}}}
\newcommand{\Ugh}{{\U}_\hb(\g)}
\newcommand{\UWurz}{\mathrm{P}}
\newcommand{\Wurz}{\Omg}

\newcommand{\z}{\mathbb Z}
\newcommand{\mm}{\mathfrak m}
\begin{document}

\title{Equivariant quantization on quotients of simple Lie groups by
reductive subgroups of maximal rank}
%
\authori{Joseph Donin and Vadim Ostapenko}
\addressi{Dept. of
Math., Bar-Ilan University, 52900 Ramat Gan, Israel}
\authorii{}     \addressii{}
\authoriii{}    \addressiii{}
\authoriv{}     \addressiv{}
\authorv{}      \addressv{}
\authorvi{}     \addressvi{}
%
\headauthor{J. Donin and V. Ostapenko} 
\headtitle{Equivariant quantization}   
\lastevenhead{J. Donin and V. Ostapenko} 
\pacs{XXX}     
\refnum{A}
\daterec{XXX}    
\issuenumber{0}  \year{2001}
\setcounter{page}{1}

\maketitle
\begin{abstract}
We consider a class of homogeneous manifolds including all 
semisimple coadjoint orbits.
We describe manifolds of that class admitting deformation q
uantizations equivariant under the action of $G$ and the 
corresponding quantum group. We also classify
Poisson brackets relating to such quantizations.
\end{abstract}
\section{Introduction}

A quantum manifold is obtained from a usual
manifold by replacing the original commutative
function algebra with a deformed non-commutative algebra.
If a manifold $M$ is equipped with an action of a Lie group $G$, 
i.e. $M$ is a $G$-manifold,
it is natural to consider deformation quantizations of the
function algebra which are equivariant with respect to the action of
the group $G$ or the corresponding quantum group $\Ugh$, 
\cite{DrinfeldV:QHopfAlgs},
where $\g$ is the Lie algebra of $G$.

We consider here homogeneous $G$-manifolds of type $M=G/K$, where
$G$ is a simple connected  Lie group over $\kk$ and $K$ its 
reductive subgroup of maximal rank, i.e. $K$ contains a maximal 
torus of $G$.
The Lie algebra of $K$, $\k$, is generated by a Cartan subalgebra of
$\g$ and the set of contained in $\k$ root vectors.
Such a class of manifolds includes semisimple coadjoint orbits of $G$ in $\g^*$.
For a coadjoint orbit, $\k$ is a Levi subalgebra generated by simple 
root vectors.

The notion of $G$-equivariant quantization is very natural and supposed to be known.
So, we only give the
definition of $\Ugh$-equivariant quantization we admit in the paper.
\begin{defn}
Let $M$ be a $G$-manifold.
A $\Ugh$-equivariant quantization on $M$
is the $\kk[[h]]$-module $\mathcal{F}(M)[[h]]$ endowed with
a deformed multiplication $\mu_h$ such that:

(i) $\mh$ defines a star-product on $M$, i.e. is of the form
$\mh=\mu+\hb\mu_1+\hb^2\mu_2+\cdots$, where $\mu$ is the usual
commutative multiplication on $\mathcal{F}(M)$ and $\mu_i$, $i\geq 1$, are bidifferential
operators on $M$ vanishing on constants;

(ii) the multiplication $\mu_h$ is $\Ugh$-equivariant, i.e.
$$
x\mh(a\otimes b)=\mh\De_\hb(x)(a\otimes b),
 \ \mbox{for any}\ x\in\Ugh,
\  a,b\in\Ah,
$$
where $\De_\hb$ denotes the comultiplication in $\Ugh$.
We also suppose that the action of $\Ugh$ on $\mathcal{F}(M)[[h]]$ is
the natural extension of the initial action of $\U(\g)$ on 
$\mathcal{F}(M)$ in the following sense.  
Since $\g$ is simple, $\Ugh$ is isomorphic to 
$\U(\g)[[h]]$ as an algebra, and
the natural action of $\U(\g)[[h]]$ on $\mathcal{F}(M)[[h]]$ is induced from
the action of $\U(\g)$ on $\mathcal{F}(M)$ by $\kk[[h]]$-linearity.
\end{defn}
The element $\nu$ defined as $\nu(a,b)=\mu_1(a,b)-\mu_1(b,a)$
for $a,b\in \mathcal{F}(M)$ is a Poisson bracket on $M$. 
In particular, $\nu$ is a bivector field on $M$. 
We call $\mu_h$ a {\em quantization corresponding} to the Poisson bracket $\nu$.
One can prove,
\cite{DoninJ:UghEquivQn},
that the Poisson bracket of any $\Ugh$-equivariant
quantization has the form
\begin{equation}\label{eqpb}
\nu=s- r_M,
\end{equation}
where $r_M$ is the bivector field which is the image of the
classical r-matrix $r\in \bigwedge^2 \g$ related to $\Ugh$ via the
action map $\g\to Vect(M)$, and $s$ is a $G$-invariant bivector
field on $M$ satisfying 
\begin{equation}\label{eq:mCYBE}
[s,s]=-\ph_M.
\end{equation}
By $\ph_M$ we denote the $3$-vector field induced on $M$
by the $\g$-invariant element $\ph=[r,r]\in\bigwedge^3\g$ and
by $[\cdot,\cdot]$ the Schouten bracket of two polyvector fields.
We call a Poisson bracket of the form (\ref{eqpb})  {\em admissible}.
\begin{rem}
It is easy to see that the admissible Poisson brackets are precisely 
those which makes $M$ a Poisson homogeneous manifold
with respect to the Drinfeld--Sklyanin Poisson bracket $r^r-r^l$ on $G$.

The admissible brackets  are also related to the dynamical CYBE
\cite{LuJH:dCYBE}, \cite{KarolinskyStolin:dCYBE}.
\end{rem}
It is convenient to introduce the following
\begin{defn}
A $G$-invariant bivector field $s$ on $M$ is called 
\emph{a $\ph$-bracket} if it satisfies $(\ref{eq:mCYBE})$.
\end{defn}

The first result of the paper  is a classification of invariant
and admissible Poisson  brackets on homogeneous manifolds under
consideration. 
Equation $(\ref{eqpb})$ reduces the classification of
admissible brackets on $M$ to a classification of
$\ph$-brackets on $M$.

The second result of the paper is that for any admissible Poisson 
bracket on $M$ there exists the corresponding $\Ugh$-equivariant
quantization.
Note that there exists a $G$- invariant connection on $M$, since the stabilizer of
a point of $M$ semisimply acts on $\g$. Moreover, since $M$ is 
$G$-homogeneous, an invariant Poisson bracket on $M$ has the same 
rank at any point of $M$.
Applying, for example, the well known Fedosov method, we obtain that for any
$G$-invariant
Poisson bracket on $M$ the corresponding $G$-equivariant 
quantization always exists.

Among other works relevant to our study, we would like to mention 
the following.
J.~Donin and  D.~Gurevich
\cite{DoninGurevic:DrJimbPoissBr} proved that on a semisimple 
coadjoint orbit there exists a quantization of
the Sklyanin--Drinfeld Poisson bracket.
J.~Donin and S.~Shnider \cite{DoninShnider:QSymSpaces}
constructed one and two parameter equivariant quantization
on any symmetric manifold.
J.~Donin, D.~Gurevich, and
S.~Shnider  \cite{DoninGurevicShnider:DoubleQuantiz}
solved the problem of one and two parameter equivariant quantization 
on any semisimple coadjoint orbit, including a complete 
classification of corresponding Poisson brackets on them.

The paper is organized as follows.
In Section~\ref{sect:RedSAlgs}, we introduce a group $\Ga$ related
to a reductive subalgebra
$\k\subset\g$ of maximal rank. It is the quotient of the lattice 
spanned by all roots of $\g$ by the sublattice spanned by roots of $\k$.
We prove a statement about how the images of roots are lying in $\Ga$.

In Section~\ref{sect:ClassifOfBrkts}, using this statement we give a 
classification of invariant Poisson brackets and $\ph$-brackets on $M$.

In Section~\ref{sect:Quantization}, we show that those brackets
can be quantized.

\section{Reductive subalgebras of maximal rank and group $\Ga$}
\label{sect:RedSAlgs}

Let $\g$ be a reductive Lie algebra over $\kk$.
Let us fix a Cartan
subalgebra $\h\subset\g$ and denote by $\Wurz\subset\hd$
the root system of $\g$.
We are interested in reductive subalgebras, $\k$, in $\g$ containing 
$\h$, which are called subalgebras of maximal rank.
To any subalgebra of this type, it corresponds a subsystem of roots,
$\UWurz$.
Namely,
$\alpha\in \UWurz$ if and only if the root vector $E_\alpha\in \k$.
Note that $\k$ is generated by
$\h$ and the root vectors $E_\alpha$, $\alpha\in\UWurz$.
\begin{eg}\label{eg:Levi}
Choose a system of simple roots $\Pi\subset\Omega$.
Fix a subset $\Sigma\subset\Pi$.
Let us put $\UWurz=\z(\Sigma)\cap\Wurz$,
where $\z(\Sigma)$ denote the free abelian group (lattice)
spanned by $\Sigma$.
Let $\k$ be the Lie subalgebra of $\g$ generated by the root
vectors $E_\al$ for all $\al\in\UWurz$ and by the Cartan
subalgebra $\h$.
Such $\k$ is called {\em a Levi subalgebra} of $\g$.
\end{eg}
\begin{eg}\label{eg:Dynkin}
Choose a system of simple roots $\Pi\subset\Omega$.
Fix a root $\al\in\Pi$ and a positive integer $n$.
Consider the set $\UWurz$ of all roots in $\Wurz$ whose coefficient
before $\al$ in their expansion in basis $\Pi$ is divisible by $n$.
Let $\k$ be the Lie subalgebra of $\g$ generated by the root
vectors $E_\al$ for all $\al\in\UWurz$ and by the Cartan
subalgebra $\h$.
\end{eg}
The following is a modified version of the well known Dynkin theorem,
\cite{DynkinE:SSSubAlgsInSSLieAlgs}.
\begin{thm}\label{DThm}
For any reductive subalgebra $\k\subset\g$ of maximal rank, 
there exists a chain 
$\g\supset\k_1\supset\k_2\supset\dots\supset\k_k=\k$ of
Lie subalgebras such that the subalgebra $\k_j$ is obtained from
$\k_{j-1}$ by applying one of the procedures described in
Examples~\ref{eg:Levi} and \ref{eg:Dynkin}.
\end{thm}
Let us set $\Ga=\z(\Wurz)/\z(\UWurz)$ and denote by $\Tor(\Ga)$ the 
subgroup of elements of finite order in $\Ga$. 
The following theorem describes the group $\Tor(\Ga)$ and images of 
roots in it.
\begin{prop}\label{MP}
\begin{description}
\item[(I)]
The group $\Ga$ is free, i.e. $\Tor(\Ga)=0$, if and only if $\k$ is 
a Levi subalgebra.
\item[(II)]
$\Tor(\Ga)$ is cyclic if and only if $\k$ is obtained as in
Theorem~\ref{DThm} applying the procedure of
Example~\ref{eg:Dynkin} no more than once.
\item[(III)]
If $\Tor(\Ga)$ is cyclic then any its non-zero element is
the image of a root.
\item[(IV)]
If $\Tor(\Ga)$ is non-cyclic, then the kernel of any character
of it contains the image of a root.
\end{description}
\end{prop}
\begin{proof}
(I) is obvious, since in this case $\UWurz$ is generated by a set of 
roots, $\UWurz'$, lying in
a linear subspace of the vector space spanned by $\Omega$. 
But it is known that such a $\UWurz'$ can be included in a set of 
simple roots.

To prove other points, we
use the structure theory of simple Lie algebras. 
First of all, we come to the following conclusions.
Any cyclic summand in $\Tor(\Ga)$ is of order $\leq 6$.
If $\g$ is a classical simple Lie algebra, the group $\Tor(\Ga)$ is
isomorphic to a direct sum of copies of $\z_2=\z/2\z$.
For exceptional simple Lie algebras $\g$, also may be 
$\Tor(\Ga)=\z_3$, $\z_4$,
$\z_5$, $\z_2\oplus\z_3$, $\z_2\oplus\z_4$, or $\z_3\oplus\z_3$.
Using this structure of $\Tor(\Ga)$,
the proof of the theorem can be completed by a direct computation.
\end{proof}

\section{Classification of invariant Poisson brackets and 
$\ph$-brackets on $M$}
\label{sect:ClassifOfBrkts}
In this section, we suppose that $G$ is a simple connected Lie group 
and $K$ a reductive subgroup of maximal rank.
We are going to describe invariant bivector fields, $s$, on $M$ 
satisfying the equation
\begin{equation}\label{eq:mCYBEkappa}
[s,s]=\ka^2\ph_M,
\end{equation}
where $\ka\in\kk$.
Note that if $s$ satisfies (\ref{eq:mCYBEkappa}) with $\ka^2=0$, 
then it defines an invariant Poisson bracket, while
if $s$ satisfies (\ref{eq:mCYBEkappa}) with $\ka^2= -1$ then it 
defines a $\ph$-bracket on $M$.

Consider the natural projection $\pi:G\to M=G/K$. It induces the map
$\pi_*:\g\to T_o$ where $T_o$ is the tangent space to $M$ at the
point $o$ being the image of unity. Since the $\ad$-action of $\k$ 
on $\g$ is semisimple,
there exists an $\ad(\k)$-invariant subspace $\mm$ of $\g$
complementary to $\k$, and one can identify $T_o$ and $\mm$ by means of
$\pi_*$. It is easy to see that subspace $\mm$ is uniquely defined
and has a basis formed by elements
$E_\gamma$, $\gamma\in\Omega\setminus P$. Here $E_\gamma$ are root
vectors satisfying $(E_\gamma,E_{-\gamma})=1$ for the Killing form $(\cdot,\cdot)$.

It follows from the above that invariant bivector fields on $M$
correspond to $\k$-invariant
bivectors of $\wedge^2\mm$.
It is proven in \cite{DoninGurevicShnider:DoubleQuantiz} and 
\cite{OstapenkoV:InvarQn}
that any bivector of $\wedge^2\mm$ is $\k$-invariant if and only if 
it has the form
\begin{equation}\label{biv}
\sum_{\alpha\in\Omega\setminus \UWurz} c(\clal)E_\alpha\wedge E_{-\alpha},
\end{equation}
where
$\clal$ is the image of $\alpha$ in $\Ga$ and
coefficients $c(\clal)$ satisfy the condition
$c(\clal)=-c(-\clal)$. This means that for invariant bivectors
the coefficients before terms $E_\alpha\wedge E_{-\alpha}$
and $E_\beta\wedge E_{-\beta}$ are the same if $\clal=\clbe$.

The following result is proved in 
\cite{DoninGurevicShnider:DoubleQuantiz}.
\begin{prop}\label{th:EqnForC}
A bivector field $s$ presented in the form (\ref{biv}) satisfies (\ref{eq:mCYBEkappa})
if and only if the
coefficients $c(\clbe)$  obey the following condition: 
if $\clal+\clbe\in\Ga$ then
$$
c(\clal+\clbe)(c(\clal)+c(\clbe))=c(\clal)c(\clbe)+\ka^2.
$$
\end{prop}
Using the previous proposition, one can obtain the following
description of $\ph$-brackets and invariant Poisson brackets on $M$,
\cite{DoninJ:QGManifolds}. See also \cite{DoninJ:UghEquivQn},
where this is done for admissible brackets on
coadjoint orbits.
\begin{prop}\label{PhiB}
A $\ph$-bracket on $M$ is determined by the following data:
a subgroup $\Psi\subset\Ga$ such that $\Ga/\Psi$ is
free, a group homomorphism
$\chi:\Psi\to\kk\setminus 0$
such that $\chi(\clal)\not=1$
for $\clal$ being the image of a root,
and a linear ordering in the group $\Ga/\Psi$.
The coefficients $c(\clal)$, $\clal\in\Ga$, from (\ref{biv}) are 
given by the formula:
$$
c(\clal)=\frac{\chi(\clal)+1}{\chi(\clal)-1} 
\qquad\mbox{for $\clal\in\Psi$,}
$$
\begin{eqnarray*}
c(\clal)&=&1 \qquad\mbox{if the projection of $\clal$ in
$\Ga/\Psi$ is positive,}\\
c(\clal)&=&-1 \qquad\mbox{if the projection of $\clal$ in
$\Ga/\Psi$ is negative.}
\end{eqnarray*}
\end{prop}
\begin{prop}\label{prop:phi}
An invariant Poisson bracket on $M$ is determined by
choosing a subgroup $\Psi\subset\Ga$ containing no torsion
elements and a group homomorphism $\la:\Psi\to\kk$ such that
$\la(\clal)\not= 0$ for  $\al$ being the image of a root.
The coefficients  $c(\clal)$ from
(\ref{biv}) are given by the formula
\begin{eqnarray*}
c(\clal)&=&\frac{1}{\la(\clal)} \qquad\mbox{for $\clal\in\Psi$,}\\
c(\clal)&=&0 \qquad\mbox{if the projection of $\clal$ in
$\Ga/\Psi$ is not zero.}
\end{eqnarray*}
\end{prop}
\begin{rem}
Another description of Poisson brackets on quotients of $G$ by
reductive subgroups of maximal rank is given in
\cite{KarolinskyStolin:dCYBE}.
\end{rem}
Combining Proposition \ref{MP}, \ref{PhiB}, and \ref{prop:phi}, 
we obtain the following
\begin{cor}\label{coro}
\begin{description}
\item[(I)]
A $\ph$-bracket (and, therefore, an admissible Poisson bracket) on $M$
exists if and only if $\Tor(\Ga)$ is cyclic.
\item[(II)]
A nonzero Poisson bracket on $M$ exists if and only if $\Ga$ is infinite.
\item[(III)]
If $\Ga$ is a finite non-cyclic group then there are neither
invariant nor admissible Poisson brackets on $M$. Therefore, such an $M$
does not admit neither $G$- nor $\Ugh$-equivariant quantization.
\end{description}
\end{cor}
\begin{proof}
(I) If $\Tor(\Ga)$ is non-cyclic then, in virtue of point (IV) of 
Proposition \ref{MP},  
$\Ker(\chi)$ must contain the image of a root, what is impossible by 
Proposition \ref{PhiB}.
On the other hand, if $\Tor(\Ga)$ is cyclic, a required homomorphism 
$\chi$ is obviously exists.

\noindent
(II) If $\Tor(\Ga)$ is finite, then $\Psi$ from Proposition 
\ref{prop:phi} must be finite too, hence
$\lambda(\Psi)=0$. It follows that $\lambda$ applying to any image 
of a root gives zero.

\noindent
(III) follows from (I) and (II).
\end{proof}
\begin{eg} 
The manifold
$M=SO(13)/(SL(2)\times SL(2)\times SL(2)\times SL(2)\times
SL(2)\times SL(2))$, has $\Ga=\z_2\times\z_2$. Therefore, it 
satisfies the point (III) of the Corollary.
\end{eg}

\section{Quantization}
\label{sect:Quantization}

\begin{thm}
Let $G$ be a semisimple Lie group, $K$ a reductive subgroup of
maximal rank. Then, for any admissible Poisson bracket on $M=G/K$ 
there exists a corresponding $\Ugh$-equivariant quantization.
\end{thm}
\begin{proof} As follows from Corollary \ref{coro}, one can suppose 
that $\Tor(\Ga)$ is cyclic. 
The case when $\Tor(\Ga)=0$, i.e. $M$ is
a semisimple coadjoint orbit, is considered in 
\cite{DoninGurevicShnider:DoubleQuantiz},
while the case when $\Ga$ is cyclic itself is considered in 
\cite{OstapenkoV:InvarQn}.
The complete proof of the theorem  obtains combining methods of the 
both papers.
\end{proof}
\begin{rem}
As was noticed in Introduction, for any invariant Poisson bracket on $M$ under consideration
there exists a corresponding $G$-equivariant quantization.
\end{rem}
\noindent
{\bf Question.}
\ Describe compatible pairs consisting of invariant and admissible
Poisson brackets on $M$.
Does there exists for any compatible pair a two parameter
$\Ugh$-equivariant quantization? This problem was solved for semisimple coadjoint
orbits in \cite{DoninGurevicShnider:DoubleQuantiz}.

\medskip
{\small {\bf Acknowledgments.\ }
We thank  A.~Elashvili, S.~Shnider, and A.~Stolin for very 
helpful discussions.}


\begin{thebibliography}{99}

\bibitem{DrinfeldV:QHopfAlgs}
V.G.Drinfeld, {\it Quasi-Hopf algebras}, Leningrad Math. J.,
1 (1990), 1419--1457.

\bibitem{DoninGurevicShnider:DoubleQuantiz}
J.Donin,  D.Gurevich, and S.Shnider,
{\it Double Quantization on Some Orbits in the Coadjoint
Representations of Simple Lie Groups},
Commun. Math. Phys., 204 (1999), no. 1, 39--60.

\bibitem{LuJH:dCYBE}
J.-H. Lu,
{\it Classical dynamical r-matrices and homogeneous Poisson structures
on $G/H$ and $K/T$}, Comm. Math. Phys., 212 (2000), no. 2, 337--370.

\bibitem{KarolinskyStolin:dCYBE}
E.Karolinsky and A.Stolin,
{\it Classical
dynamical r-matrices, Poisson homogeneous spaces,
and Lagrangian subalgebras},
math.QA/0110319.

\bibitem{DoninGurevic:DrJimbPoissBr}
J.  Donin and D.  Gurevich,
{\it Some {P}oisson structures associated to {D}rinfeld--{J}imbo
{$R$}-matrix and their quantization},
Israeli Math. Journal, 92 (1995), No 1, 23--32.

\bibitem{DoninShnider:QSymSpaces}
J.  Donin and S.  Shnider,
{\it Quantum symmetric spaces},
J. of Pure and Appl. Algebra, 100 (1995), 103--115.

\bibitem{DynkinE:SSSubAlgsInSSLieAlgs}
E.  B.  Dynkin,
{\it Semi-simple subalgebras of semi-simple {L}ie algebras},
Matemat. Sbornik N.S., 30(72) (1952),
{\sl English transl.:} {\it AMS Translations},
Series 2, vol.6  (1957), 111--244.

\bibitem{OstapenkoV:InvarQn}
V.  Ostapenko,
{\it $\Ugh$ Invariant quantization on some homogeneous manifolds},
Ph.D. Thesis (1999), Bar-Ilan University, Ramat Gan, Israel.

\bibitem{DoninJ:UghEquivQn}
J.Donin, {\it $U(g)$ equivariant quantization of coadjoint orbits
and vector bundles over them},
J. of Geometry and Physics, v. 38, 2001, 54-80.


\bibitem{DoninJ:QGManifolds}
J. Donin,
{\it Quantum $G$-manifolds},
Proc. of Workshop
"New homological and categorical methods in mathematical physics",
Manchester, July 5-12, 2001.

\end{thebibliography}
\end {document}